\documentclass[12pt,letterpaper]{article}%
\usepackage{amsmath}
\usepackage{amsfonts}
\usepackage{amssymb}
\usepackage{graphicx}%
\providecommand{\U}[1]{\protect\rule{.1in}{.1in}}
\newtheorem{theorem}{Theorem}

\newtheorem{lemma}[theorem]{Lemma}

\newtheorem{question}{Question}

\newtheorem{remark}[theorem]{Remark}

\providecommand{\boksie}{\ensuremath{\mathbin{\raisebox{0.3mm}{$\scriptstyle\square$}}}}

\begin{document}

\title{\textbf{Connected }$k$\textbf{-Dominating Graphs}}
\author{C. M. Mynhardt\thanks{Supported by the Natural Sciences and Engineering
Research Council of Canada.}{\ }and L. E. Teshima\\Department of Mathematics and Statistics\\University of Victoria, Victoria, BC, \textsc{Canada}\\{\small kieka@uvic.ca, lteshima@uvic.ca}
\and A. Roux\\Department of Mathematical Sciences\\Stellenbosch University, Stellenbosch, \textsc{South Africa}\\{\small rianaroux@sun.ac.za}}
\date{}
\maketitle

\begin{abstract}
For a graph $G=(V,E)$, the $k$-dominating graph of $G$, denoted by $D_{k}(G)$,
has vertices corresponding to the dominating sets of $G$ having cardinality at
most $k$, where two vertices of $D_{k}(G)$ are adjacent if and only if the
dominating set corresponding to one of the vertices can be obtained from the
dominating set corresponding to the second vertex by the addition or deletion
of a single vertex. We denote by $d_{0}(G)$ the smallest integer for which
$D_{k}(G)$ is connected for all $k\geq d_{0}(G)$. It is known that
$\Gamma(G)+1\leq d_{0}(G)\leq|V|$, where $\Gamma(G)$ is the upper domination
number of $G$, but constructing a graph $G$ such that $d_{0}(G)>\Gamma(G)+1$
appears to be difficult.

We present two related constructions. The first construction shows that for
each integer $k\geq3$ and each integer $r$ such that $1\leq r\leq k-1$, there
exists a graph $G_{k,r}$ such that $\Gamma(G_{k,r})=k$, $\gamma(G_{k,r})=r+1$
and $d_{0}(G_{k,r})=k+r=\Gamma(G)+\gamma(G)-1$. The second construction shows
that for each integer $k\geq3$ and each integer $r$ such that $1\leq r\leq
k-1$, there exists a graph $Q_{k,r}$ such that $\Gamma(Q_{k,r})=k$,
$\gamma(Q_{k,r})=r$ and $d_{0}(Q_{k,r})=k+r=\Gamma(G)+\gamma(G)$.

\end{abstract}


\section{Introduction}

Reconfiguration problems are concerned with determining whether -- or when --
a feasible solution to a given problem can be transformed into another such
solution via a sequence of feasible solutions in such a way that any two
consecutive solutions are adjacent according to a specified adjacency
relation. Reconfiguration versions of graph colouring and other graph
problems, such as independent sets, cliques, and vertex covers, have been
studied in e.g.~\cite{BC, CHJ1, CHJ2, Ito2, IKD}. Domination reconfiguration
problems were first considered in 2014 by Haas and Seyffarth \cite{HS1}. Their
paper stimulated the work of Suzuki, Mouawad and Nishimura \cite{SMN} as well
as their own follow-up paper \cite{HS2}.

For domination related concepts not defined here we refer the reader to
\cite{HHS}. A dominating set $D$ of a graph $G=(V,E)$ is \emph{minimal
dominating} if each vertex $v\in D$ has a $D$-\emph{private neighbour}, that
is, a vertex $v^{\prime}$ that is dominated by $v$ but by no vertex in
$D-\{v\}$. The set of $D$-private neighbours of $v$ is denoted by
$\operatorname{PN}(v,D)$. We denote the domination and upper domination
numbers of $G$ by $\gamma(G)$ and $\Gamma(G)$, respectively, and its
independence number by $\alpha(G)$. A $\gamma$-\emph{set} of $G$ is a
dominating set of $G$ of cardinality $\gamma(G)$, and a $\Gamma$-\emph{set} is
a minimal dominating set of cardinality $\Gamma(G)$. A graph $G$ is
\emph{well-dominated} if $\gamma(G)=\Gamma(G)$, \label{delete?}and
\emph{well-covered }if all its maximal independent sets have the same
cardinality $\alpha(G)$. A (not necessarily dominating) set $X\subseteq V$ is
\emph{irredundant }if $\operatorname{PN}(x,X)\neq\varnothing$ for each $x\in
X$, and we denote the maximum cardinality of an irredundant set of $G$ by
$\operatorname{IR}(G)$.

For a graph $G$, the $k$-dominating graph of $G$, denoted by $D_{k}(G)$, has
vertices corresponding to the dominating sets of $G$ having cardinality at
most $k$, and two vertices of $D_{k}(G)$ are adjacent if and only if the
dominating set corresponding to one of the vertices can be obtained from the
dominating set corresponding to the second vertex by the addition or deletion
of a single vertex. A major problem, addressed in \cite{HS1, HS2, SMN}, is to
determine conditions for $D_{k}(G)$ to be connected. As observed in
\cite{HS1}, the star $K_{1,n}$, $n\geq3$, shows that $D_{k}(G)$ being
connected does not imply that $D_{k+1}(G)$ is also connected: the unique
$\Gamma$-set of $K_{1,n}$ consists of its $n$ independent vertices and is an
isolated vertex of $D_{\Gamma}(K_{1,n})$, hence $D_{\Gamma}(K_{1,n})$ is
disconnected, but $D_{j}(K_{1,n})$ is connected for each $j$, $1\leq j\leq
n-1$. Hence an important problem in this regard is to determine the smallest
integer $j$ such that $D_{k}(G)$ is connected for all $k\geq j$.

To study this problem, let $d_{0}(G)$ denote the smallest integer for which
$D_{k}(G)$ is connected for all $k\geq d_{0}(G)$. For $G\ncong\overline{K_{n}%
}$, it is known that $\Gamma(G)+1\leq d_{0}(G)\leq|V|$, where $\Gamma(G)$ is
the upper domination number of $G$, but constructing a graph $G$ such that
$d_{0}(G)>\Gamma(G)+1$ appears to be difficult. Suzuki et al.~\cite{SMN} found
an infinite class of graphs $G_{(d,b)}$ for which $d_{0}(G_{(d,b)}%
)=\Gamma(G_{(d,b)})+2$; the smallest of these is $G_{(2,3)}\cong
P_{3}\boksie
K_{3}$ (the Cartesian product of $P_{3}$ and $K_{3}$), for which $\gamma
(P_{3}\boksie K_{3})=\Gamma(P_{3}\boksie K_{3})=3$. Haas and Seyffarth
\cite{HS2} also found a graph, which they called $G_{4}$, such that
$d_{0}(G_{4})=\Gamma(G_{4})+2$, and mentioned that they didn't know of the
existence of any graphs with $d_{0}>\Gamma+2$.

We remedy this situation by presenting two related constructions. The first
construction demonstrates the existence of graphs with arbitrary upper
domination number $\Gamma\geq3$, arbitrary domination number in the range
$2\leq\gamma\leq\Gamma$, and $d_{0}=\Gamma+\gamma-1$. (The graph
$P_{3}\boksie
K_{3}$ is an example of such a graph with $\gamma=\Gamma$.)

The second construction demonstrates the existence of graphs with arbitrary
upper domination number $\Gamma\geq3$, arbitrary domination number in the
range $1\leq\gamma\leq\Gamma-1$, and $d_{0}=\Gamma+\gamma$. For $\gamma\geq2$,
this is the first construction of graphs with $d_{0}=\Gamma+\gamma$. This
result is best possible in all cases, as it follows from \cite[Theorem 5]{HS1}
that $d_{0}(G)\leq\Gamma(G)+\gamma(G)$ for any graph $G$, and from
\cite[Theorem 7]{HS2} that $d_{0}(G)\leq2\Gamma(G)-1$ for any graph $G$. Our
main theorems are stated below.

\begin{theorem}
\label{ThmMain1}For each integer $k\geq3$ and each integer $r$ such that
$1\leq r\leq k-1$, there exists a graph $G_{k,r}$ such that $\Gamma
(G_{k,r})=k$, $\gamma(G_{k,r})=r+1$ and $d_{0}(G_{k,r})=k+r=\Gamma
(G)+\gamma(G)-1$.
\end{theorem}

\begin{theorem}
\label{ThmMain2}For each integer $k\geq3$ and each integer $r$ such that
$1\leq r\leq k-1$, there exists a graph $Q_{k,r}$ such that $\Gamma
(Q_{k,r})=k$, $\gamma(Q_{k,r})=r$ and $d_{0}(Q_{k,r})=k+r=\Gamma(G)+\gamma(G)$.
\end{theorem}

We begin Section \ref{SecThm3} by defining a new parameter $\operatorname{sep}%
(G)$ and proving in Theorem \ref{Thm3} that $\operatorname{sep}(G)=d_{0}(G)$
for all graphs $G$ with at least one edge. The constructions of $G_{k,r}$ and
$Q_{k,r}$ appear in Sections \ref{SecConstruct1} and \ref{SecConstruct2},
respectively. We then prove a number of lemmas in Section \ref{SecLemmas}; as
shown in Section \ref{SecRes} our main results are consequences of the lemmas
and Theorem \ref{Thm3}. Section \ref{SecProb} contains a list of problems for
future research.

\section{Another Definition of $d_{0}(G)$\label{SecThm3}}

A $k$\emph{-partition of a set }$A$ is a partition of $A$ into $k$ (nonempty)
subsets. Let $\mathcal{D}(G)$ be the collection of all minimal dominating sets
of $G$. For a $2$-partition $\Pi=\mathcal{X}\cup\mathcal{Y}$ of $\mathcal{D}%
(G)$, define the \emph{separation} $\operatorname{sep}(\Pi)$ \emph{in} $\Pi$
to be%
\[
\operatorname{sep}(\Pi)=\min\{|X\cup Y|:X\in\mathcal{X},\ Y\in\mathcal{Y}\}.
\]
The \emph{separation in the collection $\mathcal{D}(G)$ }of minimal dominating
sets of $G$ is
\[
\operatorname{sep}(G)=\max\{\operatorname{sep}(\Pi):\Pi\text{ is a
}2\text{-partition of }\mathcal{D}\}.
\]
If $G$ is a graph with at least one edge, then $G$ has at least two distinct
minimal dominating sets. Let $X$ be any $\Gamma$-set of $G$ and let
$\mathcal{Y}=\mathcal{D}(G)-\{X\}\neq\varnothing$. Since $\Pi^{\prime
}=\{X\}\cup\mathcal{Y}$ is a partition of $\mathcal{D}(G)$, and no minimal
dominating set has a proper subset that is a dominating set,
$\operatorname{sep}(\Pi^{\prime})>|X|=\Gamma(G)$ and thus $\operatorname{sep}%
(G)\geq\Gamma(G)+1$.

We show that $\operatorname{sep}(G)=d_{0}(G)$ for each graph $G$ having at
least one edge and use this result to prove Theorems \ref{ThmMain1} and
\ref{ThmMain2}; we expect Theorem \ref{Thm3} to be useful in other situations
as well. We state the following fact and a lemma from \cite{HS1} for referencing.

\begin{remark}
\label{RemPath}If $D$ and $D^{\prime}$ are minimal dominating sets of $G$ such
that $|D\cup D^{\prime}|\leq k$, then $D_{k}(G)$ contains a $D-D^{\prime}$
path. ($D\cup D^{\prime}$ is a vertex of $D_{k}(G)$ and there are paths from
$D\cup D^{\prime}$ to both $D$ and $D^{\prime}$.)
\end{remark}

\begin{lemma}
\label{LemHS}\cite[Lemma 4]{HS1}\hspace{0.1in}If $k>\Gamma(G)$ and $D_{k}(G)$
is connected, then $D_{k+1}(G)$ is connected.
\end{lemma}

\begin{theorem}
\label{Thm3}For any graph $G$ with at least one edge, $d_{0}%
(G)=\operatorname{sep}(G)$.
\end{theorem}

\noindent\textbf{Proof.\hspace{0.1in}}Let $k=\operatorname{sep}(G)$ and
$\mathcal{D}=\mathcal{D}(G)$. To show that $D_{k}(G)$ is connected it is
sufficient to show that $D_{k}(G)$ contains a $D-D^{\prime}$ path for any pair
of distinct minimal dominating sets $D$ and $D^{\prime}$.

Let $X_{1}$ be any minimal dominating set of $G$ and consider the partition
$\Pi_{1}=\{X_{1}\}\cup(\mathcal{D}-\{X_{1}\})$ of $\mathcal{D}$. Since
$\operatorname{sep}(\Pi_{1})\leq k$ by definition of $k$, there exists a set
$X_{2}\in\mathcal{D}-\{X_{1}\}$ such that $|X_{1}\cup X_{2}|\leq k$. By Remark
\ref{RemPath}, $D_{k}(G)$ contains an $X_{1}-X_{2}$ path. Let $\Pi_{2}%
=\{X_{1},X_{2}\}\cup(\mathcal{D}-\{X_{1},X_{2}\})$. Since $\operatorname{sep}%
(\Pi_{2})\leq k$, there exist a set $X_{i},\ i\in\{1,2\}$, and a set $X_{3}%
\in\mathcal{D}-\{X_{1},X_{2}\}$ such that $|X_{i}\cup X_{3}|\leq k$. By Remark
\ref{RemPath}, $D_{k}(G)$ contains an $X_{i}-X_{3}$ path, and so all three
sets $X_{1},X_{2},X_{3}$ belong to the same component of $D_{k}(G)$. This
process can be repeated to show that all minimal dominating sets of $G$ belong
to the same component of $D_{k}(G)$. Since any non-minimal dominating set of
cardinality at most $k$ contains a minimal dominating set and is connected to
this set in $D_{k}(G)$, $D_{k}(G)$ is connected.

Now consider $D_{k-1}(G)$; we show that it is disconnected. Since
$k>\Gamma(G)$, all minimal dominating sets of $G$ are vertices of $D_{k-1}%
(G)$. Let $\Pi=\mathcal{X}\cup\mathcal{Y}$ be a partition of $\mathcal{D}$
such that $\operatorname{sep}(\Pi)=\operatorname{sep}(G)=k$. Let
$\mathcal{X}^{\prime}$ and $\mathcal{Y}^{\prime}$ consist of all dominating
sets of $G$ of cardinality at most $k-1$ that are supersets of sets in
$\mathcal{X}$ and $\mathcal{Y}$, respectively. Then $V(D_{k-1}(G))=\mathcal{X}%
^{\prime}\cup\mathcal{Y}^{\prime}$. If there exists a dominating set
$Z\in\mathcal{X}^{\prime}\cap\mathcal{Y}^{\prime}$, then, by definition of
$\mathcal{X}^{\prime}$ and $\mathcal{Y}^{\prime}$, there also exist minimal
dominating sets $X\in\mathcal{X}$ and $Y\in\mathcal{Y}$ such that $X\subseteq
Z$ and $Y\subseteq Z$. But then $X\cup Y\subseteq Z$ and, by the choice of
$\Pi$, $k\leq|X\cup Y|\leq|Z|\leq k-1$, which is impossible. Therefore
$\mathcal{X}^{\prime}$ and $\mathcal{Y}^{\prime}$ are disjoint. Suppose
$D_{k-1}(G)$ is connected. Then some vertex $A^{\prime}\in\mathcal{X}^{\prime
}$ is adjacent to some vertex $B^{\prime}\in\mathcal{Y}^{\prime}$. By the
definition of adjacency in $D_{k-1}(G)$ we may thus assume without loss of
generality that $A^{\prime}\subset B^{\prime}$. But then there are minimal
dominating sets $A\subseteq A^{\prime}$ and $B\subseteq B^{\prime}$, implying
that $A\cup B\subseteq B^{\prime}$. But $A\in\mathcal{X}$, $B\in\mathcal{Y}$,
and $k\leq|A\cup B|\leq|B^{\prime}|\leq k-1$, which is similarly impossible.
The result now follows from Lemma \ref{LemHS}.~$\blacksquare$

\section{Constructions}

\subsection{Construction of $G_{k,r}$\label{SecConstruct1}}

\begin{figure}[ptb]%
\centering
\includegraphics[
height=2.5227in,
width=5.5884in
]%
{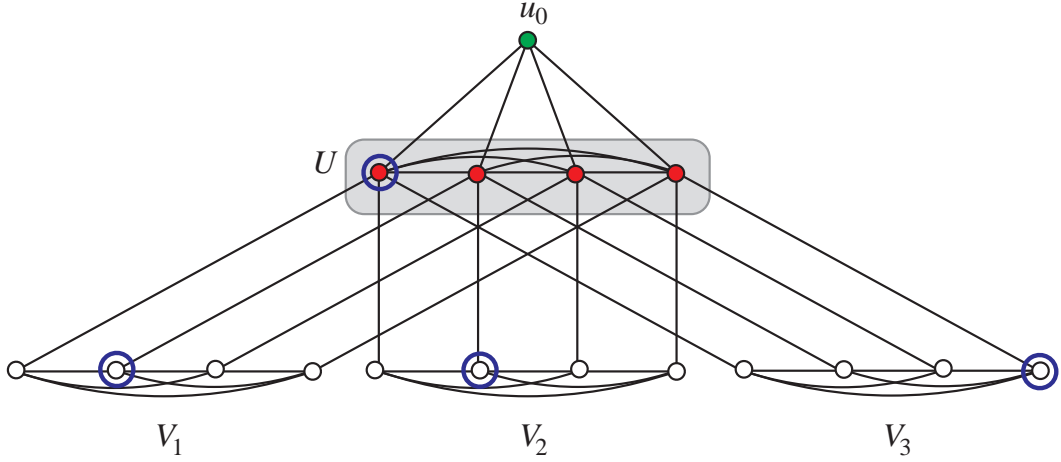}%
\caption{The well-dominated graph $G_{4,3}$ with $\gamma(G_{4,3}%
)=\Gamma(G_{4,3})=4$ and $d_{0}(G_{4,3})=7$}%
\label{Figk4r3}%
\end{figure}

(See Figure \ref{Figk4r3} for the case $k=4$ and $r=3$.) For each integer
$k\geq3$ and each integer $r$ such that $1\leq r\leq k-1$, consider
$K_{k}\boksie K_{1,r}$. If $r=1$, then $K_{1,r}=K_{2}$; in what follows
designate (rather inaccurately) one vertex of $K_{2}$ to be its centre and the
other one to be its leaf. Let $U=\{u_{1},\ldots,u_{k}\}$ (the red vertices in
Figure \ref{Figk4r3}) be the vertex set of the copy of $K_{k}$ that
corresponds to the centre of $K_{1,r}$, and for $i=1,\ldots,r$, let
$V_{i}=\{v_{i,1},\ldots,v_{i,k}\}$ be the vertex sets of the copies of $K_{k}$
that correspond to the leaves of $K_{1,r}$, where the labelling has been
chosen so that $u_{j}$ is adjacent to $v_{i,j}$ for each $i=1,\ldots,r$ and
each $j=1,\ldots,k$. Add a new vertex $u_{0}$ (the green vertex in Figure
\ref{Figk4r3}), joining it to each vertex in $U$.
This is the graph $G_{k,r}$. Let $U_{0}=U\cup\{u_{0}\}$. Define the collection
$\mathcal{X}$ of subsets of $V(G_{k,r})$ by
\[
\mathcal{X}=\{X\subseteq V(G):|X\cap U_{0}|=1\text{ and\ }|X\cap
V_{i}|=1\ \text{for each\ }i=1,\ldots,r\}.
\]
The blue circled vertices in Figure \ref{Figk4r3} form a set in $\mathcal{X}$.
Properties of dominating sets of $G_{k,r}$ are proved in Lemmas
\ref{LemDomSets}~and~\ref{CorDomSets}.

\subsection{Construction of $Q_{k,r}$\label{SecConstruct2}}

\begin{figure}[ph]%
\centering
\includegraphics[
height=2.9966in,
width=5.5582in
]%
{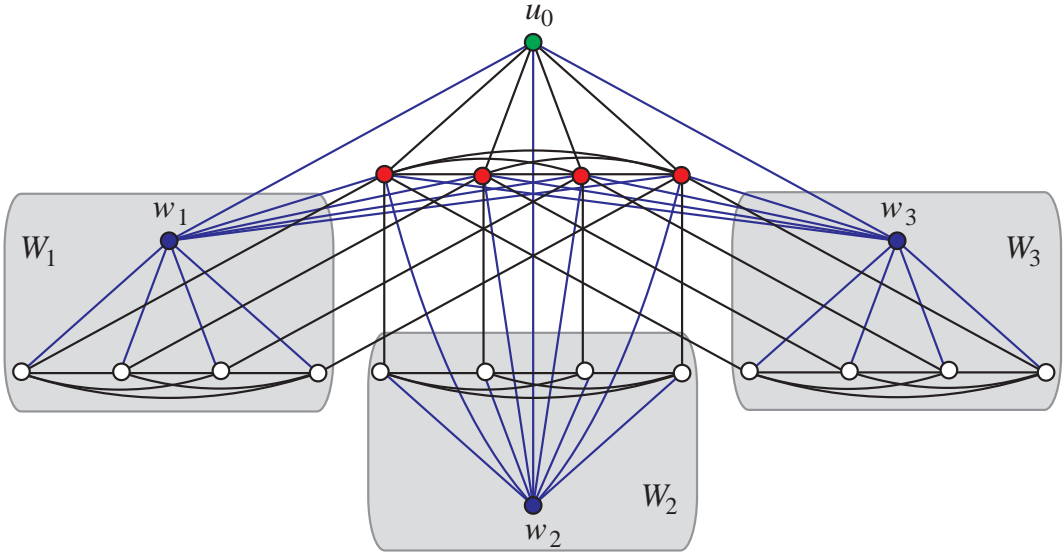}%
\caption{The graph $Q_{4,3}$ with $\gamma(Q_{4,3})=3$, $\Gamma(Q_{4,3})=4$ and
$d_{0}(Q_{4,3})=7$}%
\label{Figk4+3}%
\end{figure}

(See Figure \ref{Figk4+3} for the case $k=4$, $r=3$.) Let $G_{k,r}$ be the
graph constructed in Section \ref{SecConstruct1}, where $U$ is shown as red
vertices in Figure \ref{Figk4+3}. Add $r$ new (blue in the figure) vertices
$w_{1},\ldots,w_{r}$, joining $w_{i}$ to all vertices in $U_{0}\cup V_{i}$,
for each $i$. Call the new graph $Q_{k,r}$. Let $W_{i}=V_{i}\cup\{w_{i}\}$ and
let $\mathcal{W}$ be the collection of all subsets of $V(Q_{k,r})$ defined by%
\[
\mathcal{W}=\{Y:|Y\cap W_{i}|=1\text{ for each }i=1,\ldots,r\text{ and }%
Y\cap\{w_{1},\ldots,w_{r}\}\neq\varnothing\}.
\]
Properties of dominating sets of $G_{k,r}$ are proved in Lemmas
\ref{LemDomSetsQ}~and~\ref{LemMinDomQ}.%

\section{Lemmas\label{SecLemmas}}

In this section we prove some properties of the dominating sets of $G_{k,r}$
and $Q_{k,r}$.

\subsection{Lemmas for Theorem \ref{ThmMain1}}

\begin{lemma}
\label{LemDomSets}Let $X$ be any dominating set of $G_{k,r}$. Then

\begin{enumerate}
\item[(a)] $X\cap U_{0}\neq\varnothing$,

\item[(b)] if $U\nsubseteq X$, then $X\cap V_{i}\neq\varnothing$ for each
$i\in\{1,\ldots,r\}$.
\end{enumerate}

Assume $X$ is a minimal dominating set. Then

\begin{enumerate}
\item[(c)] $|X\cap V_{i}|\leq1$ for each $i\in\{1,\ldots,r\}$,

\item[(d)] if $u_{0}\in X$, then $X\cap U_{0}=\{u_{0}\}$,

\item[(e)] if $X\cap V_{i}=\varnothing$ for some $i$, then $X=U$.
\end{enumerate}
\end{lemma}

\noindent\textbf{Proof.\hspace{0.1in}}(a) Since $u_{0}$ is not adjacent to any
vertex in $V(G_{k,r})-U_{0}$, $X\cap U_{0}\neq\varnothing$ to dominate~$u_{0}$.

\noindent(b)\hspace{0.1in}Assume without loss of generality that $u_{k}\notin
X$. To dominate $v_{i,k}$, $X\cap V_{i}\neq\varnothing$ for each
$i\in\{1,\ldots,r\}$.

\noindent(c)\hspace{0.1in}By (a), there exists $u\in X\cap U_{0}$. If (say)
$v_{i,1},v_{i,2}\in X$, then $v_{i,1}$ and $v_{i,2}$ have no private
neighbours in $V_{i}$. Hence their private neighbours are in $U$. This is
impossible because all vertices in $U$ are dominated by $u$.

\noindent(d)\hspace{0.1in}If $u_{0}\in X$, then $\operatorname{PN}(u_{0},X)\in
U_{0}$. Since $U_{0}$ is a clique, $U\cap X=\varnothing$.

\noindent(e)\hspace{0.1in}If $X\cap V_{i}=\varnothing$ for some $i$, then by
(b), $U\subseteq X$. Since $U$ dominates $\bigcup_{j=1}^{r}V_{j}$, no vertex
$z\in\bigcup_{j=1}^{r}V_{j}$ has a $U\cup\{z\}$-private neighbour, hence
$X\cap(\bigcup_{j=1}^{r}V_{j})=\varnothing$. By (d), $u_{0}\notin X$. It
follows that $X=U$.~$\blacksquare$

\begin{lemma}
\label{CorDomSets}If $r<k-1$, then $U$ is the only $\Gamma$-set of $G_{k,r}$,
and $\mathcal{X}$ is precisely the collection of $\gamma$-sets of $G$. If
$r=k-1$, then $G$ is well-dominated and the above-mentioned sets are precisely
the $\gamma$-sets (and $\Gamma$-sets) of $G$. Moreover, $G$ has no other
minimal dominating sets.
\end{lemma}

\noindent\textbf{Proof.\hspace{0.1in}}Since $U$ is a minimal dominating set of
$G_{k,r}$, $\Gamma(G_{k,r})\geq k$. Suppose $X\neq U$ is a minimal dominating
set of $G_{k,r}$. By Lemma \ref{LemDomSets}(e), $X\cap V_{i}\neq\varnothing$
for each $i\in\{1,\ldots,r\}$, and by Lemma \ref{LemDomSets}(c), $|X\cap
V_{i}|=1$ for each $i$. By Lemma \ref{LemDomSets}(a), $X\cap U_{0}%
\neq\varnothing$. Since $X\cap V_{i}\neq\varnothing$ for each $i\in
\{1,\ldots,r\}$ and $G_{k,r}[V_{i}]$, the subgraph of $G_{k,r}$ induced by
$V_{i}$, is complete, $\operatorname{PN}(u,X)\subseteq U_{0}$ for each $u\in
X\cap U_{0}$. Since $G_{k,r}[U_{0}]=K_{k+1}$, $|X\cap U_{0}|=1$. Moreover,
$u_{0}\in\operatorname{PN}(u,X)$ for the unique vertex $u\in X\cap U_{0}$.
Hence $\mathcal{D}(G_{k,r})=\mathcal{X}\cup\{U\}$. The result follows because
$|U|=k$ and $|X|=r+1\leq k$ for each $X\in\mathcal{X}$.~$\blacksquare$

\bigskip

\noindent\textbf{Note:\hspace{0.1in}}The set $\{u_{1},\ldots,u_{k-1}%
\}\cup\{v_{1,k},v_{2,k},\ldots,v_{r-1,k}\}$ is a non-dominating irredundant
set of $G_{k,r}$, hence $\operatorname{IR}(G_{k,r})\geq k-1+r-1$.

\subsection{Lemmas for Theorem \ref{ThmMain2}}

Lemma \ref{LemDomSets} can be adapted to $Q_{k,r}$ as follows.

\begin{lemma}
\label{LemDomSetsQ}Let $X$ be any dominating set of $Q_{k,r}$. Then

\begin{enumerate}
\item[(a)] $X\ $contains at least one vertex in $U_{0}\cup\{w_{1},\ldots
,w_{r}\}$.

\item[(b)] If $U\nsubseteq X$, then $X\cap W_{i}\neq\varnothing$ for each
$i\in\{1,\ldots,r\}$. For each $i\in\{1,\ldots,r\}$, if $U\nsubseteq X$ and
$w_{i}\notin X$, then $X\cap V_{i}\neq\varnothing$.
\end{enumerate}

Assume $X$ is a minimal dominating set. Then

\begin{enumerate}
\item[(c)] $|X\cap W_{i}|\leq1$ for each $i\in\{1,\ldots,r\}$.

\item[(d)] If $u_{0}\in X$, then $X\cap U_{0}=\{u_{0}\}$.

\item[(e)] If $X\cap W_{i}=\varnothing$ for some $i$, then $X=U$.

\item[(f)] If $w_{i}\in X$ for some $i$, then $X\cap U_{0}=\varnothing$.
\end{enumerate}
\end{lemma}

\noindent\textbf{Proof.\hspace{0.1in}}We only prove (f) as the proofs of (a)
-- (e) are similar to the corresponding proofs of Lemma \ref{LemDomSets}. Note
that $U$ and each $Y\in\mathcal{W}$ are (minimal) dominating sets of $Q_{k,r}%
$. Say $w_{1}\in X$. By the minimality of $X$, $U\nsubseteq X$, hence by (b),
$X\cap W_{i}\neq\varnothing$ for each $i\in\{1,\ldots,r\}$. But then
$Y\subseteq X$ for some $Y\in\mathcal{W}$ and, again by the minimality of $X$,
$X=Y$ and so $X\cap U_{0}=\varnothing$.~$\blacksquare$

\begin{lemma}
\label{LemMinDomQ}The minimal dominating sets of $Q_{k,r}$ are precisely the
sets in $\mathcal{W}$, which are the $\gamma$-sets, the sets in $\mathcal{X}$,
which are $\Gamma$-sets if and only if $r=k-1$, and $U$, the unique $\Gamma
$-set if $r<k-1$.
\end{lemma}

\noindent\textbf{Proof.\hspace{0.1in}}It is easy to check that any
$Y\in\mathcal{W}$ dominates $Q_{k,r}$, hence $\gamma(Q_{k,r})\leq r$. Let $Z$
be any minimal dominating set of $Q_{k,r}$. If $Z\cap W_{i}=\varnothing$ for
some $i$, then $Z=U$ by Lemma \ref{LemDomSetsQ}(e). Hence assume $Z\cap
W_{i}\neq\varnothing$ for each $i$.

Suppose $u,u^{\prime}\in Z\cap U_{0}$. Since $Q_{k,r}[U_{0}]\ $is complete,
the $Z$-private neighbours of $u$ (and those of $u^{\prime}$) belong to
$W_{i}$ for some $i$. But this is impossible because each $Q_{k,r}[W_{i}]$ is
also complete and $Z\cap W_{i}\neq\varnothing$ for each $i$. Therefore $|Z\cap
U_{0}|\leq1$.

If $Z\cap\{w_{1},\ldots,w_{r}\}=\varnothing$, then $Z\cap U_{0}\neq
\varnothing$ in order to dominate $u_{0}$. By Lemma \ref{LemDomSetsQ}(c) and
the results above, $|Z\cap U_{0}|=1$ and $|Z\cap V_{i}|=1$ for each $i$, that
is, $Z\in\mathcal{X}$. On the other hand, if $Z\cap W_{i}\neq\varnothing$ for
each $i$ and $Z\cap\{w_{1},\ldots,w_{r}\}\neq\varnothing$, then $Z\cap
U_{0}=\varnothing$ by Lemma \ref{LemDomSetsQ}(f), hence $Z\in\mathcal{W}$.
Therefore $\mathcal{D}(Q_{k,r})=\mathcal{X}\cup\mathcal{W}\cup\{U\}$.

Since $|Y|=r$ for each $y\in\mathcal{W}$, $|X|=r+1$ for each $X\in\mathcal{X}%
$, and $r<k=|U|$, it follows that $\gamma(Q_{k,r})=r$, $\Gamma(Q_{k,r})=k$,
the $\gamma$-sets are precisely the sets in $\mathcal{W}$, and the $\Gamma
$-sets consist of $U$ and, if $r=k-1$, also the sets in $\mathcal{X}%
$.~$\blacksquare$

\section{Proofs of Main Results\label{SecRes}}

Our main results, which we reformulate here for convenience, now follow easily.

\noindent\textbf{Theorem \ref{ThmMain1}\hspace{0.1in}}\emph{For each integer
}$k\geq3$\emph{ and each integer }$r$\emph{ such that }$1\leq r\leq
k-1$\emph{, the graph }$G_{k,r}$\emph{ satisfies }$\Gamma(G_{k,r})=k$\emph{,
}$\gamma(G_{k,r})=r+1$\emph{ and }$d_{0}(G_{k,r})=k+r$\emph{.}

\bigskip

\noindent\textbf{Proof.\hspace{0.1in}}By Lemma \ref{CorDomSets},
$\mathcal{D}(G_{k,r})=\mathcal{X}\cup\{U\}$.\textbf{ }Consider the
$2$-partition $\Pi=\mathcal{X}\cup\{U\}$ of $\mathcal{D}(G_{k,r})$. For
$X\in\mathcal{X}$, if $u_{0}\in X$, then $|X\cup U|=k+r+1$, otherwise $|X\cup
U|=k+r$. Hence $\operatorname{sep}(\Pi)=k+r$ and $\operatorname{sep}%
(G_{k,r})\geq k+r$. Consider any other $2$-partition $\Pi^{\prime}%
=\mathcal{X}^{\prime}\cup\mathcal{Y}^{\prime}$ of $\mathcal{D}(G_{k,r})$ where
$U\in\mathcal{Y}^{\prime}$, say. If there exists $X^{\prime}\in\mathcal{X}%
^{\prime}$ such that $u_{0}\notin X^{\prime}$, then $|X^{\prime}\cup U|=k+r$
and $\operatorname{sep}(\Pi^{\prime})\leq k+r$. Hence assume $u_{0}$ belongs
to each set in $\mathcal{X}^{\prime}$. For any $X^{\prime}\in\mathcal{X}%
^{\prime}$, let $Z=(X^{\prime}-\{u_{0}\})\cup\{u_{1}\}$ and note that
$Z\in\mathcal{Y}^{\prime}$. Since $|X^{\prime}\cup Z|=r+2\leq k+1$, it follows
that $\operatorname{sep}(\Pi^{\prime})\leq k+1\leq k+r$. This shows that
$\operatorname{sep}(G_{k,r})=k+r$ and the result follows from Theorem
\ref{Thm3}.~$\blacksquare$

\bigskip

Haas and Seyffarth \cite{HS2} show that $d_{0}=\Gamma+1$ for large classes of
well-covered and well-dominated graphs. The graphs $G_{(d,b)}$ in \cite{SMN}
are well-covered with $d_{0}(G_{(d,b)})=\Gamma(G_{(d,b)})+2$, while $G_{4}$ in
\cite{HS2} is well-dominated with $d_{0}(G_{4})=\Gamma(G_{4})+2$. However, it
is not true that well-covered and well-dominated graphs have relatively small
values of $d_{0}$: consider any $\gamma$-set $X\in\mathcal{X}$ of $G_{k,r}$
such that $X\cap U_{0}=\{u_{0}\}$. Then $X$ is an independent set; indeed,
since $X$ is dominating, it is a maximal independent set. Since $U$ is not
independent, and $U$ and the sets in $\mathcal{X}$ are the only minimal
dominating sets of $G_{k,r}$ (by Lemma \ref{CorDomSets}), $\alpha
(G_{k,r})=\gamma(G_{k,r})=r+1$, and $G_{k,r}$ is well-covered. Therefore
$d_{0}(G_{k,r})=\Gamma(G)+\alpha(G)-1$. Hence Theorem \ref{ThmMain1} also
shows that the bound in the first part of Theorem 7 of \cite{HS2} is tight,
even for well-covered graphs. If $r=k-1$, then $G_{k,r}$ is well-dominated,
and presents an example of a well-dominated graph such that $d_{0}%
(G_{k,r})=2\Gamma(G)-1$, thus equalling the upper bound for $d_{0}$.

Since $\alpha(G_{k,r})=\gamma(G_{k,r})$, $d_{0}(G_{k,r})=\Gamma(G_{k,r}%
)+\gamma(G_{k,r})-1$. It is therefore tempting to think that $d_{0}%
(G)\leq\Gamma(G)+\gamma(G)-1$ for all graphs $G$ with $\gamma(G)\geq2$, but
Theorem \ref{ThmMain2} shows that this is not the case.

\bigskip

\noindent\textbf{Theorem \ref{ThmMain2}\hspace{0.1in}}\emph{For each integer
}$k\geq3$\emph{ and each integer }$r$\emph{ such that }$1\leq r\leq
k-1$\emph{, the graph }$Q_{k,r}$\emph{ satisfies }$\Gamma(Q_{k,r})=k$\emph{,
}$\gamma(Q_{k,r})=r$\emph{ and }$d_{0}(Q_{k,r})=k+r$\emph{.}

\bigskip

\noindent\textbf{Proof.\hspace{0.1in}}By Lemma~\ref{LemMinDomQ},
$\mathcal{D}(Q_{k,r})=\mathcal{X\cup W}\cup\{U\}$. Consider the $2$-partition
$\Pi=(\mathcal{X\cup W)}\cup\{U\}$ of $\mathcal{D}(Q_{k,r})$. For any
$X\in\mathcal{X}$, $k+r\leq|X\cup U|\leq k+r+1$, as in the proof of Theorem
\ref{ThmMain1}. For any $Y\in\mathcal{W}$, $|Y\cup U|=k+r$. Hence
$\operatorname{sep}(\Pi)=k+r$, so that $\operatorname{sep}(Q_{k,r})\geq k+r$.
By Theorem \ref{Thm3}, $d_{0}(Q_{k,r})\geq k+r$. Since $d_{0}(G)\leq
\gamma(G)+\Gamma(G)$ for all graphs $G$ \cite[Theorem 5]{HS1}, the result
follows.~$\blacksquare$

\section{Future Research\label{SecProb}}

All graphs $G$ such that $d_{0}(G)>\Gamma(G)+1$ constructed so far contain
triangles (in fact, their clique numbers increase with increased upper
domination number). Since $\alpha=\Gamma$ for all bipartite graphs, Theorem 9
in \cite{HS2} shows that $d_{0}(G)=\Gamma(G)+1$ if $G$ is bipartite. The
problem remains open for triangle-free non-bipartite graphs.

\begin{question}
Is it true that $d_{0}(G)=\Gamma(G)+1$ if $G$ is triangle-free?
\end{question}

Suzuki et al.~\cite{SMN} constructed an infinite family of graphs $G_{n}$ of
order $63n-6$ such that $D_{\gamma(G)+1}(G_{n})$ has exponential diameter
$\Omega(2^{n})$. They then posed the question of whether there exists a value
of $k$ for which $D_{k}(G)$ is connected and guaranteed not to have
exponential diameter (in the order of $G$). Since $D_{|V(G)|}(G)$ is connected
and has diameter at most $2(|V(G)|-\gamma(G))$, such a $k$ certainly exists,
and we pose the following questions instead.

\begin{question}

\begin{enumerate}
\item[$(i)$] If $G$ is an $n$-vertex graph and $k=d_{0}(G)$, what is the
maximum/minimum diameter of $D_{k}(G)$? Which $n$-vertex graphs attain this maximum/minimum?

\item[$(ii)$] Suppose $D_{i}(G)$ and $D_{j}(G)$ are connected and $i<j$. How
are $\operatorname{diam}(D_{i}(G))$ and $\operatorname{diam}(D_{j}(G))$
related? (If $i>\Gamma(G)$, then $\operatorname{diam}(D_{i}(G))\geq
\operatorname{diam}(D_{j}(G))$.)
\end{enumerate}
\end{question}

The following problems are mentioned in \cite{HS1}. We repeat them here
because they are interesting and worthy of consideration.

\begin{question}
\cite{HS1}\hspace{0.1in}When is $D_{k}(G)$ Hamiltonian?
\end{question}

Note that $D_{k}(G)$ is bipartite for any graph $G$ and any integer
$k\geq\gamma(G)$, since even-sized dominating sets are only adjacent to
odd-sized sets and vice versa. In fact, since dominating sets are subsets of
$V(G)$, and the definition of adjacency is the same in $D_{k}(G)$ and in the
hypercube $Q_{n}$ whose vertices represent the subsets of an $n$-set,
$D_{k}(G)$ is an induced subgraph of $Q_{|(V(G)|}$. For example, $D_{n}%
(K_{n})$ is the hypercube $Q_{n}$ with the vertex corresponding to the empty
set deleted, while $D_{2}(K_{1,n})\cong K_{1,n}$, $n\geq3$. These observations
lead to our next questions.

\begin{question}

\begin{enumerate}
\item[$(i)$] (Adapted from \cite{HS1}) Which induced subgraphs of $Q_{n}$
occur as $D_{k}(G)$ for some $n$-vertex graph $G$ and some integer $k$?

\item[$(ii)$] \cite{HS1} For which graphs $G$ is $D_{k}(G)\cong G$ for some
value of $k$?
\end{enumerate}
\end{question}

\begin{question}
\cite{HS1}\hspace{0.1in}What is the complexity of determining whether two
dominating sets of $G$ are in the same component of $D_{k}(G)$, or of
$D_{\Gamma(G)+1}(G)$?
\end{question}

\end{document}